\documentclass[12pt]{article}

\usepackage{amsmath, epsfig, cite}
\usepackage{amssymb}
\usepackage{amsfonts}
\usepackage{latexsym}
\usepackage{float}
\usepackage{color}

\newtheorem{thm}{Theorem}[section]

\newtheorem{lem}[thm]{Lemma}

\numberwithin{equation}{section}

\newcommand{\qed}{{\hfill$\square$}\medskip}

\begin{document}

\begin{center}
{\Large\bf Supercongruences for Almkvist--Zudilin sequences}
\end{center}

\vskip 2mm \centerline{Ji-Cai Liu$^1$ and He-Xia Ni$^2$}
\begin{center}
{\footnotesize $^1$Department of Mathematics, Wenzhou University, Wenzhou 325035, PR China\\
{\tt jcliu2016@gmail.com } \\[10pt]

$^2$Department of Applied Mathematics, Nanjing Audit University, Nanjing 211815, PR China\\
{\tt nihexia@yeah.net  }
}
\end{center}


\vskip 0.7cm \noindent{\bf Abstract.}
In this note, we prove two supercongruences involving Almkvist--Zudilin sequences, which were originally conjectured by Z.-H. Sun.


\vskip 2mm
\noindent{\it MR Subject Classifications}: 11A07, 05A19, 11B68

\section{Introduction}
In 2006, Almkvist and Zudilin \cite{az-b-2006} introduced the following Ap\'ery-like sequences:
\begin{align*}
G_n=\sum_{k=0}^n{2k\choose k}^2{2n-2k\choose n-k}4^{n-k}.
\end{align*}
Recently, Z.-H. Sun \cite{sunzh-a-2020} investigated congruence properties of the above sequences.
For instance, Z.-H. Sun \cite[Theorems 2.2 and 2.6]{sunzh-a-2020} showed that for any prime $p\ge 5$,
\begin{align}
&G_{p-1}\equiv (-1)^{\frac{p-1}{2}}256^{p-1}\label{a-1}\\[10pt]
&+p^2\left(E_{p-3}-8(-1)^{\frac{p-1}{2}}q_p(2)^2+\frac{1}{2}
\sum_{k=0}^{(p-1)/2}\frac{{2k\choose k}^2}{16^k}H_k^2\right)\pmod{p^3},\notag\\[10pt]
&\sum_{k=0}^{p-1}\frac{G_k}{16^k}\equiv p^2\left(1-\sum_{k=0}^{(p-1)/2}\frac{{2k\choose k}^2}{16^k(k+1)}H_k\right)\pmod{p^3}.\label{a-2}
\end{align}
Here $q_p(2)$ is the Fermat quotient $(2^{p-1}-1)/p$, the $n$th harmonic number is given by
\begin{align*}
H_n=\sum_{k=1}^n\frac{1}{k},
\end{align*}
and the Euler numbers are defined as
\begin{align*}
\frac{2}{e^x+e^{-x}}=\sum_{n=0}^{\infty}E_n\frac{x^n}{n!}.
\end{align*}

The motivation of this note is to refine \eqref{a-1} and \eqref{a-2} by establishing the following results, which were originally conjectured by Z.-H. Sun \cite[Conjectures 2.1 and 2.2]{sunzh-a-2020}.
\begin{thm}\label{t-1}
For any prime $p\ge 5$, we have
\begin{align}
&G_{p-1}\equiv (-1)^{\frac{p-1}{2}}256^{p-1}+3p^2E_{p-3}\pmod{p^3},\label{a-4}\\[10pt]
&\sum_{k=0}^{p-1}\frac{G_k}{16^k}\equiv p^2\left(4(-1)^{\frac{p-1}{2}}-3\right)\pmod{p^3}.\label{a-5}
\end{align}
\end{thm}

\section{Proof of Theorem \ref{t-1}}
In order to prove Theorem \ref{t-1}, we require the following identities.
\begin{lem}
For any non-negative integer $n$, we have
\begin{align}
&\sum_{k=0}^{n}(-1)^k{n\choose k}{n+k\choose k}H_k^2=
2(-1)^n\left(2H_n^2+\sum_{k=1}^n\frac{(-1)^k}{k^2}\right),\label{b-1}\\[10pt]
&\sum_{k=0}^{n}\frac{(-1)^k}{k+1}{n\choose k}{n+k\choose k}H_k=\frac{(-1)^n-1}{n(n+1)}.\label{c-1}
\end{align}
\end{lem}
{\noindent \it Proof.}
Identities \eqref{b-1} and \eqref{c-1} can be discovered and proved by the symbolic summation
package {\tt Sigma} developed by Schneider \cite{schneider-slc-2007}. One can refer to \cite{liu-jsc-2019} for the same approach to finding and proving identities of this type.
In particular, \eqref{b-1} was also proved by Wang \cite[Lemma 2.2]{wang-jmaa-2020}.
\qed

{\noindent\it Proof of \eqref{a-4}.}
Letting $n=\frac{p-1}{2}$ in \eqref{b-1} and noting that
\begin{align}
&(-1)^k{\frac{p-1}{2}\choose k}{\frac{p-1}{2}+k\choose k}\equiv \frac{{2k\choose k}^2}{16^k}\pmod{p},\label{new-1}
\end{align}
we obtain
\begin{align}
\sum_{k=0}^{(p-1)/2}\frac{{2k\choose k}^2}{16^k}H_k^2\equiv
2(-1)^{\frac{p-1}{2}}\left(2H_{\frac{p-1}{2}}^2+\sum_{k=1}^{(p-1)/2}\frac{(-1)^k}{k^2}\right)\pmod{p}.\label{b-2}
\end{align}

By \cite[(30)]{lehmer-am-1938} and \cite[Lemma 2.4]{sunzw-scm-2011}, we have
\begin{align}
&H_{\frac{p-1}{2}}\equiv -2q_p(2)\pmod{p},\label{b-3}\\[10pt]
&\sum_{k=1}^{(p-1)/2}\frac{(-1)^k}{k^2}\equiv 2(-1)^{\frac{p-1}{2}}E_{p-3}\pmod{p}.\label{b-4}
\end{align}
It follows from \eqref{b-2}--\eqref{b-4} that
\begin{align}
\sum_{k=0}^{(p-1)/2}\frac{{2k\choose k}^2}{16^k}H_k^2\equiv
16(-1)^{\frac{p-1}{2}}q_p(2)^2+4E_{p-3}\pmod{p}.\label{b-5}
\end{align}
Substituting \eqref{b-5} into \eqref{a-1}, we arrive at \eqref{a-4}.
\qed

{\noindent \it Proof of \eqref{a-5}.} Letting $n=\frac{p-1}{2}$ in \eqref{c-1} and using
\eqref{new-1}, we obtain
\begin{align}
\sum_{k=0}^{(p-1)/2}\frac{{2k\choose k}^2}{16^k(k+1)}H_k\equiv 4\left(1-(-1)^{\frac{p-1}{2}}\right)\pmod{p}.\label{c-3}
\end{align}
Substituting \eqref{c-3} into \eqref{a-2}, we complete the proof of \eqref{a-5}.
\qed

\vskip 5mm \noindent{\bf Acknowledgments.}
The authors are grateful to Professor Zhi-Wei Sun for his valuable comments on a previous version of this paper. The first author was supported by the National Natural Science Foundation of China (grant 11801417).

\end{document}